\documentclass[12pt]{amsart}

\usepackage[utf8]{inputenc}
\usepackage{amsmath}
\usepackage{amsfonts,amssymb,latexsym,multicol}
\usepackage[francais]{babel}
\usepackage[T1]{fontenc} 
\usepackage{amsthm}

\newcounter{fig}
\def\figcaption #1
  {\addtocounter{fig}{1}
   \begin{center}
   {\sc Fig.} \arabic{fig}. {\it #1}
    \end{center}
    \addtocontents{lof}{\protect \contentsline {figure}{\protect \numberline {\thefig}{\protect \ignorespaces \protect \pit #1 }}{\thepage}}
   }
\input{cyracc.def}

\newtheorem{theo}{Th\'eor\`eme}
\newtheorem{lem}{Lemme}

\newtheorem{prop}{Proposition}

\newcommand{\expli}[1]{\quad\text{\footnotesize (#1)}}

\newcommand{\ioe}{\leqslant}
\newcommand{\soe}{\geqslant}
\newcommand{\vers}{\rightarrow}

\newcommand{\demi}{{\frac{1}{2}}}

\newcommand{\Ccal}{{\mathcal C}}
\newcommand{\Dcal}{{\mathcal D}}

\newcommand{\Lcal}{{\mathcal L}}

\newcommand{\cgot}{{\mathfrak c}}

\newcommand{\Nat}{{\mathbb N}}
\newcommand{\Int}{{\mathbb Z}}
\newcommand{\Rat}{{\mathbb Q}}
\newcommand{\Real}{{\mathbb R}}
\newcommand{\Com}{{\mathbb C}}

\newcommand{\sgn}{{\rm sgn}}

\newcommand{\fin}{\hfill$\Box$}
\newcommand{\dem}{\noindent {\bf D\'emonstration\ }}
\let\oldqedhere\qedhere
    \renewcommand{\qedhere}{\pushQED{\qed}\oldqedhere}

\title{Comportement local moyen de la fonction de Brjuno}
\author{Michel Balazard et Bruno Martin}
\thanks{Bruno Martin a b\'en\'efici\'e d'un financement par la fondation autrichienne
FWF dans le cadre du projet S9605 de son r\'eseau Analytic Combinatorics and
Probabilistic Number Theory. Les auteurs remercient d'autre part le laboratoire
Poncelet (CNRS-Universit\'e Ind\'ependante de Moscou) pour son accueil lors de la
pr\'eparation de cet article.}
\begin{document}
\maketitle


\begin{center}
  {\sc Abstract}
\end{center}
\begin{quote}
{\footnotesize We describe the average behaviour of the Brjuno function $\Phi$ in the neighbourhood of any given point of the unit interval. In particular, we show that the Lebesgue set of $\Phi$ is the set of Brjuno numbers and we find the asymptotic behaviour of the modulus of continuity of the integral of $\Phi$.}
\end{quote}

\begin{center}
  {\sc Keywords}
\end{center}
\begin{quote}
{\footnotesize Brjuno function, Brjuno numbers, Lebesgue set, Continued fractions. }
\end{quote}

\section{Introduction}\label{intro}

Désignons par $\lfloor u \rfloor$ la partie enti\`ere du nombre r\'eel $u$, et par $\{u\}$ sa partie fractionnaire :
$$
\{u\}=u-\lfloor u \rfloor, \quad \lfloor u \rfloor \in \Int, \quad 0\ioe \{u\} <1.
$$
Posons  $X=]0,1[ \setminus \Rat$. L'application
\begin{align*}
\alpha : X &\vers X \\
x & \mapsto \{1/x\}
\end{align*}
définit une transformation mesurable de $X$ muni de la tribu bor\'elienne dans lui-m\^eme. Ryll-Nardzewski (cf. \cite{MR0046583}, Theorem 2, p.76) a observ\'e que la mesure de probabilit\'e
$$
\mu=\frac{dx}{(1+x)\log 2}
$$
est invariante\footnote{On trouve une assertion \'equivalente dans l'article de Doeblin \cite{MR0002732}, p.357.} par $\alpha$ : la mesure image $\alpha\mu$ n'est autre que $\mu$.

En d'autres termes, si $f \in L^1(0,1)$, alors $f \circ \alpha$  appartient aussi \`a $L^1(0,1)$ et
$$
\int_0^1 f\bigl (\alpha(x)\bigr )\frac{dx}{1+x}=\int_0^1f(x)\frac{dx}{1+x}.
$$
Nous définissons les fonctions it\'er\'ees de $\alpha$ en posant $\alpha_0(x)=x$ et pour $k\soe 1$, $\alpha_k(x)=\alpha \big(\alpha_{k-1}(x)\big)$.
Les fonctions $\alpha_k\,(k\soe0)$ sont définies sur $X$, \`a valeurs dans $X$, et laissent la mesure $\mu$ invariante\footnote{On peut montrer que $\alpha$ est une transformation \emph{m\'elangeante} sur $(X,\mu)$. Cela permet de démontrer une assertion de Gauss dans une lettre \`a Laplace du 30 janvier 1812 (sur la convergence de la fonction de r\'epartition de $\alpha_k$, par rapport \`a la mesure de Lebesgue, quand $k \vers \infty$, cf. \cite{46.1426.04}, p.371-374). Cette assertion de Gauss avait \'et\'e démontr\'ee pour la premi\`ere fois par Kuzmin en 1928, cf.\cite{Kuzmin}. Le \S 4 de \cite{MR0192027} expose clairement cette th\'eorie.}.
\medskip

Posons
\begin{equation}\label{defPhi}
\Phi(x)=\sum_{k\soe 0}\alpha_0(x)\alpha_1(x)\cdots\alpha_{k-1}(x)\log\big(1/\alpha_k(x)\big).
\end{equation}
Cette s\'erie \`a termes positifs définit une fonction $\Phi :  X \vers ]0,\infty]$, appel\'ee fonction de Brjuno\footnote{Il s'agit ici d'une variante de la fonction définie par Yoccoz dans \cite{Yoccoz}, pp. 11-16 : nous utilisons la th\'eorie habituelle des fractions continues, alors que Yoccoz en considére une version l\'eg\'erement modifi\'ee.}. On constate ais\'ement que $\Phi$ satisfait \`a l'\'equation fonctionnelle
\begin{equation}\label{eq_fonctionnelle}
\Phi(x)=\log(1/x)+x\Phi\big(\alpha(x)\big) \quad(x\in X).
\end{equation}
Les nombres irrationnels $x$ pour lesquels $\Phi(x)=\infty$, appel\'es nombres de Cremer, constituent un ensemble de mesure de Lebesgue nulle.
On appelle nombres de Brjuno les nombres $x$ pour lesquels la s\'erie converge.

La fonction de Brjuno intervient dans l'\'etude des syst\'emes dynamiques engendr\'es par les it\'erations d'une fonction holomorphe $f:\Com\mapsto\Com$. Citons en particulier un th\'eor\`eme de Yoccoz \cite{Yoccoz} : si $f'(0)=e^{2i\pi x}$ o\`u $x\in X$, alors $f$ est localement lin\'earisable au voisinage de $z=0$ si, et seulement si $x$ est un nombre de Brjuno.

Marmi, Moussa et Yoccoz \cite{Marmi_Moussa_Yoccoz} ont entrepris l'\'etude de la r\'egularit\'e de la fonction de Brjuno et de certaines de ses variantes. Ils \'etablissent notamment que $\Phi$ appartient \`a $L^p(0,1)$ pour tout $p\soe 1$ et, ce qui est une propri\`et\'e plus forte, qu'elle est \`a oscillation moyenne born\'ee :
\begin{equation}
 \sup \frac{1}{|I|}\int_I \Big| \Phi(x) -\frac{1}{|I|}\int_I \Phi(t)dt\Big|dx <\infty,
\end{equation}
o\`u le supremum est pris sur tous les intervalles ouverts $I$ inclus dans $]0,1[$, et o\`u $|I|$ d\'esigne la longueur de $I$.

\bigskip

Dans ce travail, nous caract\'erisons les points de Lebesgue de $\Phi$ : \'etant donn\'ee une fonction $f$ localement int\'egrable sur un intervalle $J$ de $\Real$, on dit que $x$ est un point de Lebesgue de $f$ si
\begin{equation*}
  \frac{1}{h}\int_{x-h/2}^{x+h/2} \rvert f(t)-f(x) \lvert dt =o(1)\quad(h\vers
0).
\end{equation*}
Un th\'eor\`eme classique de la th\'eorie de l'int\'egrale de Lebesgue stipule que
presque tout \'el\'ement de $J$ est un point de Lebesgue de $f$ (cf.
\cite{wheeden-zygmund} Chapitre 7 par exemple). Notre r\'esultat principal est le
suivant.
\begin{theo}\label{thprincipal}
  Les points de Lebesgue de la fonction $\Phi$ sont exactement les nom\-bres de Brjuno.
\end{theo}

Notre deuxi\`eme r\'esultat est une description compl\`ete du comportement local de la fonction croissante et absolument continue $\Psi$ définie par
\begin{equation}\label{defPsi}
\Psi(x)=\int_0^x \Phi(t)dt \quad (0 \ioe x \ioe 1).
\end{equation}

Pour le formuler, nous prolongeons d'abord $\Phi$ aux nombres rationnels de
l'intervalle $]0,1[$. Si $r \in ]0,1[\,\cap\, \Rat$, on peut définir
$\alpha_k(r)$ par it\'eration de $\alpha$ tant que $k$ ne dépasse pas la longueur
$K=K(r)$ de l'\'ecriture de $r$ en fraction continue ; nous appellerons profondeur
de $r$ cette quantit\'e $K(r)$ (voir au \S \ref{profondeur_rationnel} ci-dessous
la définition pr\'ecise). On a alors $\alpha_K(r)=0$ et $\alpha_k(r)>0$ pour
$0\ioe k<K$. On pose
$$
\Phi(r) = \sum_{0\ioe k<K}\alpha_0(r)\alpha_1(r)\cdots\alpha_{k-1}(r)\log\big(1/\alpha_k(r)\big).
$$
On a par exemple $\Phi(1/k)=\log k$ si $k$ est entier, $k\soe 2$. Par convention, on pose $\Phi(0)=\Phi(1)=0$. L'\'equation fonctionnelle \eqref{eq_fonctionnelle} est alors encore valable si $x\in ]0,1[\,\cap\, \Rat$. 

\begin{theo}\label{t68}
 Soit $x_0 \in [0,1]$. On a
\begin{align*}
(i) \quad & \Psi(x_0+h)-\Psi(x_0)=o\bigl(h \log(1/|h|)\bigr) \quad (h \vers 0,
\, x_0 \in X)\\
(ii) \quad & \frac{\Psi(x_0+h)-\Psi(x_0)}{h}\vers \Phi(x_0)\quad (h \vers 0, \, x_0 \in X)\\
(iii) \quad & \Psi(x_0+h)-\Psi(x_0)\\&=\frac{1}{q} h\log (1/|h|)
+\Big( \frac{1}{q}-2\frac{\log q}{q}+\Phi(x_0)\Big)h+O\bigl(qh^2 \log
((q^2|h|)^{-1})\bigr) \\& (x_0=\frac{p}{q}, \, p \in \Nat, \, q \in \Nat^*,
 (p,q)=1, \, |h|<2/3q^2),
\end{align*}
o\`u la constante implicite dans le $O$ est absolue.
\end{theo}

\medskip

Soit
$$
\omega(h) = \sup \{|\Psi(x)-\Psi(y)|, \, 0 \ioe x,y \ioe 1, \, |x-y| \ioe h\}
$$
le module de continuit\'e de la fonction $\Psi$. Notre troisi\`eme r\'esultat fournit un \'equivalent asymptotique pour $\omega (h)$ quand $h$ tend vers $0$.
\begin{theo}\label{t69}
On a
$$
\omega(h) = h \log(1/h) +O(h) \quad (0<h\ioe 1).
$$
\end{theo}

\bigskip

Nous emploierons désormais les notations
\begin{equation}\label{t63}
 \beta_k(x)=\alpha_0(x)\alpha_1(x)\cdots\alpha_{k}(x)\end{equation}
(par convention, $\beta_{-1}=1$) et
 \begin{equation}\label{defgamma}
 \gamma_k(x)=\beta_{k-1}(x)\log\big(1/\alpha_k(x)\big)\quad(x\in X, k\soe 0),
 \end{equation}
de sorte que $\gamma_0(x)=\log (1/x)$ et
que
\begin{align}
\Phi(x)&=\sum_{k\soe 0} \gamma_k(x)\\
&= \sum_{k < K} \gamma_k(x)+\beta_{K-1}(x)\Phi\big(\alpha_K(x)\big)
 \quad(x\in X, \, K \in \Nat).\label{t61}
\end{align}
L'identit\'e \eqref{t61} est une forme g\'en\'eralis\'ee de l'\'equation fonctionnelle \eqref{eq_fonctionnelle}.

Le \S \ref{t80} contient quelques remarques g\'en\'erales sur le comportement local moyen des fonctions localement int\'egrables. Au \S \ref{pfractions_continues}, nous rappelons les \'el\'ements de la th\'eorie des fractions continues qui nous seront utiles.  Nous obtenons au \S \ref{t81} une premi\`ere majoration du module de continuit\'e $\omega(h)$, essentiellement gr\^ace au th\'eor\`eme de Ryll-Nardzewski. Pour la suite, notre approche est fondée sur un découpage diophantien de $X$ qui fait intervenir la notion
de cellule, déj\`a sous-jacente dans \cite{Marmi_Moussa_Yoccoz}. Nous définissons
cette notion et l'\'etudions en détail au \S \ref{cellule}, en introduisant
notamment les notions de profondeur et d'\'epaisseur d'un sous-intervalle de
$]0,1[$. Nous \'evaluons ensuite l'int\'egrale des fonctions $\gamma_k$ sur un
intervalle au \S \ref{pmajorations_moyenne}.  Au \S \ref{t77} nous d\'emontrons
le point $(iii)$ du th\'eor\`eme \ref{t68}. Apr\`es avoir trait\'e le cas des
nombres de Cremer au \S \ref{etude_Psi} nous d\'emontrons le th\'eor\`eme
\ref{thprincipal} au \S \ref{pth_principal} et le th\'eor\`eme \ref{t69} au \S
\ref{module-continuite}. Au \S \ref{t87} final, nous formulons des questions
ouvertes.

Toutes les constantes implicites dans les symboles $O$ sont absolues. Cela \'etant, les calculs conduisent naturellement \`a des majorations num\'eriquement explicites, et nous avons généralement adopté ce mode de présentation. Enfin, nous emploierons la notation d'Iverson : $[A]=1$ si la propri\'et\'e $A$ est v\'erifi\'ee, $[A]=0$ sinon.

\section{G\'en\'eralit\'es}\label{t80}

Soit $I$ un intervalle de $\Real$, $X$ une partie de $I$ telle que $I\setminus X$ est de mesure nulle, et $f : X \vers \Real$ une fonction int\'egrable sur chaque segment inclus dans $I$. Pour fixer les id\'ees, nous supposerons que $0 \in I$.

Consid\'erons les quatre ensembles suivants :

$\bullet$ l'ensemble $\Ccal$ des points de $X$ o\`u $f$ est continue ;

$\bullet$ l'ensemble $\Lcal$ des points $x_0 \in X$ tels que
$$
\frac{1}{h}\int_{x_0-h/2}^{x_0+h/2} \rvert f(x)-f(x_0) \lvert dx \vers
0\quad(h\vers 0) \, ;
$$

$\bullet$ l'ensemble $\Lcal^1$ des points $x_0 \in I$ tels qu'il existe $\alpha \in \Real$ avec
$$
\frac{1}{h}\int_{x_0-h/2}^{x_0+h/2} \rvert f(x)-\alpha \lvert dx \vers
0\quad(h\vers 0) \, ;
$$

$\bullet$ l'ensemble $\Dcal$ des points de d\'erivabilit\'e de la fonction
$$
F(x) =\int_0^xf(t)dt \quad (x \in I).
$$

L'ensemble $\Lcal$ est l'ensemble de Lebesgue de $f$, et $\Lcal^1$ est l'{\og ensemble de Lebesgue $L^1$\fg} de $f$ (voir par exemple \cite{MR1133575}, \S 2). On a
$$
\Ccal \subset \Lcal \subset \Lcal^1 \subset \Dcal,
$$
chaque inclusion pouvant \^etre stricte. L'ensemble $I\setminus \Lcal$ est de mesure nulle (th\'eor\`eme de Lebesgue).

Le cas de la fonction de Brjuno, que nous analysons en d\'etail dans cet article, illustre la question g\'en\'erale suivante. Soit $f_n : I \vers [0, \infty[$ ($n \soe 0$) une suite de fonctions int\'egrables telle que
$$
f(x)=\sum_{n \soe 0}f_n(x)
$$
soit int\'egrable sur $I$. Comment d\'eterminer les ensembles $\Ccal$, $\Lcal$, $\Lcal^1$ et $\Dcal$ \`a partir de la connaissance des fonctions $f_n$?

Sans entrer plus avant dans l'\'etude de ce probl\`eme en toute
g\'en\'eralit\'e, con\-ten\-tons-nous d'indiquer un exemple montrant qu'un point
$x_0 \in I$ tel que :

$\bullet$ chaque fonction $f_n$ est continue en $x_0$\\
et

$\bullet$ la s\'erie $f(x_0)$ converge,\\
n'est pas n\'ecessairement un point de d\'erivabilit\'e de $F$. Il suffit de prendre $I=\Real$, $f_0=0$ et, pour $n \soe 1$,
$$
f_n(t)=|t|\cdot [|t|\ioe n^{-2/3}].
$$

Chaque $f_n$ est continue en $0$ avec $f_n(0)=0$, donc $f(0)=0$. Cependant $f(x)=|x|^{-1/2}-|x|\{|x|^{-3/2}\}$ ($x\not =0$), donc $F$ n'est pas d\'erivable en $0$.

Si $f$ est la fonction de Brjuno $\Phi$, on a $\Ccal =\emptyset$ (la fonction de Brjuno est non born\'ee au voisinage de tout point de $[0,1]$, au vu, par exemple, du point $(iii)$ du th\'eor\`eme \ref{t68}) et $\Lcal=\Lcal^1=\Dcal$ est l'ensemble des nombres de Brjuno, d'apr\`es le th\'eor\`eme \ref{thprincipal} et le point $(ii)$ du th\'eor\`eme \ref{t68}.

\section{Fractions continues}\label{pfractions_continues}

Dans cette section nous rappelons quelques \'el\'ements de la th\'eorie des fractions continues, et nous donnons des estimations \'el\'ementaires concernant les fonctions $\alpha_k$ et $\beta_k$.

\subsection{Quotients incomplets, r\'eduites}
 Pour $k \in \Nat$, nous définissons la fonction\footnote{Pour tout ce qui concerne ce paragraphe on pourra se r\'ef\'erer par exemple \`a \cite{Niven_Montgomery}, chapter 7.}
$$
a_k : X \vers \Nat
$$
par les relations $a_0 = 0$ et $a_k= \lfloor 1/\alpha_{k-1} \rfloor \quad (k \soe 1)$.
Les fonctions $a_k$, $k \soe 1$, sont \`a valeurs dans $\Nat^*$. On appelle $a_k(x)$ le {\og $k$\up{e} quotient incomplet\fg} de $x$.
Avec la notation classique
$$
[X_0;X_1,\dots,X_k]=X_0 +\cfrac{1}{X_1+
                          \cfrac{1}{X_2+
                           \cfrac{1}{\ddots  +
                            \cfrac{1}{X_{k-1}+
                            \cfrac{1}{X_k
                        }}}}}
,$$
on a pour tout $k \soe 0$,
$$
x=[a_0(x);a_1(x),\dots,a_{k-1}(x),a_k(x)+\alpha_k(x)] \quad (x \in X).
$$

On définit aussi les fonctions
$
p_k, \, q_k : X \vers \Nat
$
par les relations
\begin{align*}
  p_{-1}&=1 &q_{-1}&=0\\
p_0 &=0 & q_0&=1\\
p_1&=1& q_1&=a_1\\
&\vdots &&\vdots\\
p_k&=a_kp_{k-1}+p_{k-2}& q_k&=a_kq_{k-1}+q_{k-2} \quad (k \soe 1)
\end{align*}
de sorte que
$$
[a_0;a_1,\dots,a_{k-1},t]=\frac{p_{k-1}t+p_{k-2}}{q_{k-1}t+q_{k-2}} \quad (k\soe 1, \, t >0),
$$
et
$$
[a_0;a_1,\dots,a_k]=\frac{p_k}{q_k}.
$$
Les fonctions $q_k$ ($k \soe 0$) et $p_k$ ($k \soe 1$) sont \`a valeurs dans $\Nat^*$. La fraction $p_k(x)/q_k(x)$ est appel\'ee la {\og $k$\up{e} r\'eduite\fg} de $x$. Elle est irr\'eductible car
\begin{equation}\label{pkqk1}
p_kq_{k-1}-q_kp_{k-1}=(-1)^{k-1}.
\end{equation}
Une autre identit\'e dont nous aurons l'usage est
\begin{equation}\label{pkqk2}
p_{k+1}q_{k-1} -p_{k-1} q_{k+1}=(-1)^{k-1}a_{k+1}.
\end{equation}
La suite $\{q_k\}_{k\soe -1}$ est \`a croissance au moins
exponentielle\footnote{Nous entendons par l\`a l'existence d'une constante $c>1$
telle que $q_k\soe c^k$. Mais on n'a pas forc\'ement $\liminf_{k\vers
\infty}q_{k+1}/q_k >1$.}. En effet, on a pour tout $k \soe 1$
\begin{equation}\label{croissance_qk}
q_k=a_kq_{k-1}+q_{k-2} \soe q_{k-1}+ q_{k-2}\soe 2q_{k-2}
\end{equation}
En particulier,
\begin{equation}
  \label{mino_qk}
  q_k \soe F_{k+1}  \quad (k\soe 0),
\end{equation}
o\`u $F_n$ est le $n$\up{e} nombre de Fibonacci ($F_0=0$, $F_1=1$ et $F_{n+1}=F_{n-1}+F_n$ pour $n \soe 1$).
On déduit \'egalement de \eqref{croissance_qk} que
\begin{align}
  \sum_{j=0}^kq_j&= \sum_{j=2}^kq_{j-2}+q_{k-1}+q_k\notag\\
&\ioe \sum_{j=2}^k(q_j-q_{j-1})+q_{k-1}+q_k\notag\\
&=2q_k+q_{k-1}-q_1\notag\\
&\ioe 3q_k,\label{somme_qk}
\end{align}
majoration valable pour tout $k \in \Nat$.

\smallskip

Nous utiliserons \`a plusieurs reprises une majoration de la somme des inverses des nombres de Fibonacci :
\begin{equation}\label{t82}
\sum_{k\soe 0} 1/F_ {k+1}\ioe 3,36.
\end{equation}

\bigskip

\subsection{Fonctions $\alpha_k$ et $\beta_k$}

De la formule \begin{align}
  x &= [a_0(x);a_1(x),\dots,a_{k-1}(x),a_k(x)+\alpha_k(x)] \notag \\
&= \frac{p_k(x)+\alpha_k(x)p_{k-1}(x)}{q_k(x)+\alpha_k(x)q_{k-1}(x)} \label{t64},
\end{align}
on déduit
\begin{equation}
  \label{formule_alphak}
\alpha_k(x)=-\frac{p_k(x)-xq_k(x)}{p_{k-1}(x)-xq_{k-1}(x)}.
\end{equation}
Cela implique les identit\'es
\begin{align}
\beta_k(x)&=(-1)^{k-1}\bigl(p_k(x)-xq_k(x)\bigr )\label{formule_beta}\\
&= \frac{1}{q_{k+1}(x)+\alpha_{k+1}(x)q_{k}(x)} \quad(x\in X).\label{t62}
\end{align}
On dispose par cons\'equent de l'encadrement
\begin{equation}\label{encabeta}
 \frac{1}{q_{k+1}+q_k} \ioe \beta_k \ioe \frac{1}{q_{k+1}} \quad (k \soe -1).
\end{equation}
Une majoration plus g\'en\'erale nous sera utile : la définition \eqref{t63} entra\^ine
\begin{equation}\label{betajk}
 \beta_{i+j}=\beta_i\cdot \beta_{j-1}\circ\alpha_{i+1}\quad(i \soe -1, j \soe 0),
\end{equation}
d'ou l'on déduit, d'apr\`es \eqref{encabeta} et \eqref{mino_qk},
\begin{equation}\label{majobetakj}
\beta_{i+j} \ioe \frac{1}{q_{i+1}F_{j+1}}.
\end{equation}
Par ailleurs, notant que $\lfloor 1/\alpha_k\rfloor=a_{k+1}$ et que $a_{k+1}=(q_{k+1}-q_{k-1})/q_k$, nous avons
d'apr\`es \eqref{croissance_qk}
\begin{equation}\label{encaalpha}
 \frac{q_{k+1}}{2q_k}\ioe \frac{1}{\alpha_{k}} \ioe a_{k+1}+1 \quad (k \soe 0).
\end{equation}

\subsection{Comparaison des fonctions $\gamma_k$ et $q_k^{-1}\log q_{k+1}$}  

Les fonctions $\gamma_k$ et $q_k^{-1}\log q_{k+1}$ sont proches l'une de l'autre. La proposition suivante exprime pr\'ecis\'ement ce fait.
\begin{prop}\label{t0}
Pour $k\soe 0$, on a
\begin{equation}\label{enca_gamma}
-\frac{\log (2q_k)}{q_k}\ioe \gamma_k-\frac{\log q_{k+1}}{q_k}\ioe
[k=0]\frac{\log
2}{q_k}. 
\end{equation}
\end{prop}
\dem

D'apr\`es \eqref{t62}, on a 
\[
\frac{\log q_{k+1}}{q_k} -\gamma_k = \frac{q_k \log(q_{k+1}\alpha_k)
+\alpha_kq_{k-1}\log q_{k+1} }{q_k(q_k+\alpha_kq_{k-1}}
\soe \frac{q_k \log (q_{k+1} \alpha_k)}{q_k(q_k+\alpha_kq_{k-1})}.
\]
Or pour $k>0$, $q_{k-1}\soe 1$, et donc
\[
 q_{k+1}\alpha_k= \frac{ a_{k+1} q_k+q_{k-1}}{a_{k+1}+\alpha_{k+1}}
\soe  1,
\]
On obtient ainsi la deuxi\`eme in\'egalit\'e  de \eqref{enca_gamma} lorsque
$k>0$, c'est-\'e-dire, 
\[
\gamma_k\ioe \frac{\log q_{k+1}}{q_k}.
\]
Pour $k=0$, on a
\[ \gamma_0 -\frac{\log q_1}{q_0}=
  \log(1/x) -\log a_1
= \log\bigl( \tfrac{a_1+\alpha_1}{a_1}\bigr)
\ioe \log\bigl( \tfrac{a_1+1}{a_1}\bigr) \ioe \log 2.
\]

Par ailleurs, 
$$
\gamma_k=\frac{\log (1/\alpha_k)}{q_k+\alpha_kq_{k-1}}
$$
est une fonction décroissante de $\alpha_k$. Comme $\alpha_k\ioe 1/a_{k+1}$, on a donc
$$
\gamma_k\soe \frac{\log a_{k+1}}{q_k+q_{k-1}/a_{k+1}}=\frac{a_{k+1}\log a_{k+1}}{q_{k+1}}.
$$
Pour montrer que
$$
\gamma_k-\frac{\log q_{k+1}}{q_k}\soe -\frac{\log (2q_k)}{q_k},
$$
il suffit donc de montrer que
\begin{equation}
  \label{t1}
  a_{k+1}\log a_{k+1}\soe \frac{q_{k+1}}{q_k}\log \frac{q_{k+1}}{2q_k}.
\end{equation}

Si $a_{k+1}=1$, on a $q_{k+1}=q_k+q_{k-1}\ioe 2q_k$, donc \eqref{t1} est v\'erifi\`ee. Si $a_{k+1}\soe 2$, on a $\frac{q_{k+1}}{q_k}\soe 2$ et
$$
a_{k+1}=\frac{q_{k+1}-q_{k-1}}{q_k}\soe \frac{q_{k+1}}{q_k} -1.
$$
L'in\'egalit\'e \eqref{t1} r\'esultera alors de
$$
(x-1)\log (x-1)\soe x\log (x/2) \quad (x\soe 2)
$$
ce qui découle de l'\'etude du sens de variation de la diff\'erence
\[(x-1)\log
(x-1)- x\log (x/2)\] sur $[2,\infty[$.\fin

\medskip

En particulier, $x\in X$ est un nombre de Brjuno si et seulement si la s\'erie
\begin{equation}
 \sum_{k\soe0} \frac{\log q_{k+1}(x)}{q_k(x)}
\end{equation}
est convergente.

\section{Majoration du module de continuit\'e de $\Psi$}\label{t81}

Sans utiliser la notion de cellule (\S \ref{cellule}), on peut déj\`a donner une majoration utile du module de continuit\'e de $\Psi$, proposition \ref{t67} ci-dessous. Pour la démontrer, nous majorons d'abord l'int\'egrale de $\log   (1/\alpha_k) $ sur un intervalle de longueur $h$, proposition \ref{majo_log_alphak} ci-dessous.

Notre démonstration de la proposition \ref{majo_log_alphak}
requiert une majoration de la fonc\-tion $\Gamma$ d'Euler, sans doute classique mais pour laquelle nous ne disposons pas de r\'e\-f\'e\-ren\-ce. Rappelons que  pour $x>0$,
\begin{equation*}
  \Gamma(x)=\int_0^\infty e^{-t} t^{x-1} d t=\int_0^1(\log 1/t)^{x-1}dt.
\end{equation*}
\begin{lem}\label{majo_Gamma}
Pour $q\in[2,\infty[$, on a
$$
2\Gamma(q+1)\ioe q^q.
$$
\end{lem}
\dem

Notons $\psi$ la dériv\'ee logarithmique de la fonction $\Gamma$.
L'in\'egalit\'e propos\'ee est vraie pour $q=2$. Ensuite, il suffit de voir que les dériv\'ees logarithmiques des deux membres v\'erifient la m\'eme in\'egalit\'e pour $q \soe 2$, c'est-\'e-dire que
\begin{equation}\label{majo_psi}
\psi (q+1) \ioe \log  q +1 \quad (q \soe 2).
\end{equation}

Or, nous disposons de la formule classique
 \begin{equation}
\psi (x) = -\gamma-\frac{1}{x} +\sum_{n \soe 1}\frac{1}{n}-\frac{1}{n+x} \quad (x \not = 0, -1, -2, \dots),
 \end{equation}
o\`u $\gamma$ désigne la constante d'Euler.
Pour $x>0$, on a donc
\begin{align*}
\psi (x) &= -\gamma-\frac{1}{x} +\sum_{n \soe 1}\frac{x}{n(n+x)}\\
&\ioe -\gamma-\frac{1}{x} +1-\frac{1}{x+1} + \int_1^{\infty} \Bigl (\frac{1}{t}-\frac{1}{t+x} \Bigr ) dt\\
&= 1-\gamma-\frac{1}{x} -\frac{1}{x+1} + \log (x+1).
 \end{align*}
Par cons\'equent,
\begin{align*}
\psi (q+1)-\log q &\ioe \log \Bigl ( 1+\frac{2}{q}\Bigr ) -\frac{1}{q+1}-\frac{1}{q+2} +1 -\gamma\\
&\ioe 1 \quad (q \soe 2).\qedhere
 \end{align*}

\begin{prop}\label{majo_log_alphak}
Soit $I$ un intervalle de longueur $h \ioe e^{-2}$, inclus dans $[0,1]$. On a pour tout $k \in \Nat$
$$
\int_I \log  \bigl (1/\alpha_k(x)\bigr ) dx\ioe eh\log (1/h).
$$
\end{prop}
\dem

Pour $q>1$ et $p=q/(q-1)$, l'in\'egalit\'e de H\"older donne
\begin{align*}
\int_I \log  \bigl (1/\alpha_k(x)\bigr ) dx \ioe h^{1/p} \Bigl ( \int_0^1  \log ^q \bigl (1/\alpha_k(x)\bigr )  dx \Bigr )^{1/q}.
\end{align*}

Comme $\alpha$ laisse la mesure $\mu$ invariante, nous avons
\begin{align*}
\int_0^1  \log ^q \bigl (1/\alpha_k(x)\bigr )  dx
&\ioe 2 \int_0^1 \log ^q \bigl (1/\alpha_k(x)\bigr ) \frac{dx}{1+x} = 2 \int_0^1 \log ^q (1/x) \frac{dx}{1+x}\\
& \ioe 2\int_0^1 \log ^q (1/x) dx = 2 \Gamma(q+1).
\end{align*}
D'apr\`es le lemme \ref{majo_Gamma}, il suit
\begin{equation*}
\int_I \log  \bigl (1/\alpha_k(x)\bigr ) dx \ioe
qh^{1-1/q} \quad(q\soe 2).
\end{equation*}
En effectuant le choix $q=\log(1/h)$, pour lequel $q\soe 2$ d'apr\`es l'hypoth\'ese faite sur $h$, nous obtenons la majoration annonc\'ee.
\fin

\begin{prop}\label{t67}
On a
$$
\omega (h) \ioe 10 h\log (1/h) \quad (0<h \ioe e^{-2}).
$$
\end{prop}
\dem

On a pour $0<h \ioe e^{-2}$, 
\begin{align*}
 \Psi(x+h)-\Psi(x)&= \int_x^{x+h} \Phi(t) dt =
\sum_{k\ge 0} \int_x^{x+h} \gamma_k(t) dt \\
&\ioe \sum_{k\ge 0} \frac{1}{F_{k+1}}\int_x^{x+h}
\log(1/\alpha_k(t))dt \expli{d'apr\`es \eqref{encabeta} puis \eqref{mino_qk}}\\
&\ioe eh \log(1/h) \sum_{k\ge 0} \frac{1}{F_{k+1}}\expli{d'apr\`es la proposition
\ref{majo_log_alphak}}\\
&\ioe 3,36e h\log(1/h) \expli{d'apr\`es \eqref{t82}},
\end{align*}
et cela entra\^ine bien la majoration annonc\'ee. \fin 

\section{Notion de cellule}\label{cellule}

\subsection{Définition et propri\`et\'es}

Soit $b_0=0$, $b_1, \dots, b_k \in\Nat^*$. La cellule (de profondeur $k$) $\cgot(b_1,\dots,b_k)$ est l'intervalle ouvert d'extr\'emit\'es $[b_0;b_1,\dots,b_k]$ et $[b_0;b_1,\dots,b_{k-1},b_k+1]$. Par convention, la seule cellule de profondeur $0$ est $]0,1[$.

La notion de cellule n'est pas nouvelle : elle intervient naturellement dans le cadre de la th\'eorie m\'etrique et ergodique des fractions continues sous l'appelation de cylindre (cf \cite{Dajani} page 30 et chapitre 3.5) ou d'intervalle (fondamental) d'ordre $n$ (cf \cite{MR1451873} \S 12 ou \cite{MR0192027} \S 4 par exemple).

Nous \'enon\c{c}ons \`a pr\'esent quelques propri\'et\'es \'el\'ementaires concernant les cellules.

$\bullet$ Dans la cellule $\cgot(b_1,\dots,b_k)$, les fonctions $a_j$, $p_j$, $q_j$ sont constantes  pour $j \ioe k$ :
$$
a_j(x)=b_j, \quad \frac{p_j(x)}{q_j(x)}=[b_0;b_1,\dots,b_j] \quad (x \in \cgot(b_1,\dots,b_j)).
$$
\smallskip

$\bullet$ La cellule $\cgot(b_1,\dots,b_k)$ est l'intervalle ouvert d'extr\'emit\'es
$$
\frac{p_k}{q_k} \quad \text{et} \quad \frac{p_k+p_{k-1}}{q_k+q_{k-1}}
$$
(dans cet ordre si $k$ est pair ; dans l'ordre oppos\'e si $k$ est impair). Sa longueur est
\begin{equation}
  \label{t83}
\frac{1}{q_k(q_k+q_{k-1})}.
\end{equation}

\smallskip

$\bullet$ On a
$$
\alpha \bigl (\cgot(b_1,\dots,b_k)\cap X\bigr )=\cgot(b_2,\dots,b_k)\cap X.
$$

\smallskip

$\bullet$ Soit $k\soe 1$, et $\cgot= \cgot(b_1 \dots, b_{k})$ une cellule de profondeur $k$.
Les cellules de profondeur $k+1$ incluses dans $\cgot$ forment une suite $(\cgot_n)_{n \soe 1}$ d'intervalles ouverts qui constitue une partition de $\cgot$ priv\'ee d'une partie de $\Rat$ (les rationnels de la forme $[b_0;b_1,\dots,b_{k},n]$, $n$ décrivant $\Nat^*$).

La cellule $\cgot_n$ est l'intervalle d'extr\'emit\'es
\begin{align*}
 [b_0;b_1,\dots,b_k,n] &= \frac{np_k+p_{k-1}}{nq_k+q_{k-1}}\\
\text{et} \quad [b_0;b_1,\dots,b_{k-1},b_k,n+1] &=\frac{(n+1)p_k+p_{k-1}}{(n+1)q_k+q_{k-1}}.
\end{align*}
Ainsi $\cgot_n$ et $\cgot_{n+1}$ sont contig\"ues, et $\cgot_n$ est l'ensemble des nombres de la forme
$$
\frac{sp_k+p_{k-1}}{sq_k+q_{k-1}}, \quad n<s<n+1.
$$
Dans $\cgot_n$, la fonction $q_{k+1}$ a la valeur constante $nq_k+q_{k-1}$.
Observons que $\frac{p_k+p_{k-1}}{q_k+q_{k-1}}$ est l'une des extr\'emit\'es de $\cgot_1$, mais que $\frac{p_k}{q_k}$ n'est l'extr\'emit\'e d'aucune cellule $\cgot_n$. En revanche, $\frac{p_k}{q_k}$ est la limite quand $n$ tend vers l'infini des extr\'emit\'es de $\cgot_n$.

\subsection{Distance d'un irrationnel donn\'e au bord de la cellule de profondeur $k$ qui le contient}

Soit $x \in X$ et $k \in \Nat$. Il existe une unique cellule de profondeur $k$ qui contient $x$ :
$$
\cgot = \begin{cases}
\displaystyle{\Big]\frac{p_k}{q_k} ,\frac{p_k+p_{k-1}}{q_k+q_{k-1}}\Big[} \quad \text{si $k$ est pair}\\ \\
\Big]\displaystyle{\frac{p_k+p_{k-1}}{q_k+q_{k-1}},\frac{p_k}{q_k}\Big[} \quad \text{si $k$ est impair,}
\end{cases}
$$
o\`u l'on a pos\'e $p_j=p_j(x)$, $q_j=q_j(x)$ pour $j \soe -1$.
La distance de $x$ au bord de $\cgot$ joue un r\^ole déterminant dans la suite de notre travail. Nous la notons $\delta_k(x)$. On a par exemple
$$
\delta_0(x)=\min(x,1-x) \quad (x \in X).
$$
Observons tout de suite que la suite $k \mapsto \delta_k(x)$ est décroissante (au sens large) et tend vers $0$ quand $k$ tend vers l'infini. Cherchons maintenant \`a exprimer $\delta_k(x)$ en termes des fonctions pr\'ec\'edemment définies. On a d'abord,
d'apr\`es \eqref{formule_beta},
\begin{equation}\label{distance1}
 x-\frac{p_k}{q_k} = (-1)^k\frac{\beta_k(x)}{q_k}.
\end{equation}
Pour la distance de $x$ \`a l'autre extr\'emit\'e de $\cgot$, on a
\begin{align*}
x-\frac{p_k+p_{k-1}}{q_k+q_{k-1}}&=(-1)^{k-1}\frac{\beta_{k+1}(x)}{q_{k+1}}+\Bigl (\frac{p_{k+1}}{q_{k+1}}-\frac{p_k+p_{k-1}}{q_k+q_{k-1}}\Big)\\
&=(-1)^{k-1} \Bigl (  \frac{\beta_{k+1}(x)}{q_{k+1}} + \frac{a_{k+1}-1}{q_{k+1}(q_k+q_{k-1})} \Bigr ),
\end{align*}
o\`u la deuxi\`eme \'egalit\'e r\'esulte des identit\'es \eqref{pkqk1} et \eqref{pkqk2}. Nous obtenons ainsi la formule
$$
\delta_k= \min \Bigl ( \frac{\beta_k}{q_k}, \frac{\beta_{k+1}}{q_{k+1}} + \frac{a_{k+1}-1}{q_{k+1}(q_k+q_{k-1})} \Bigr ).
$$
\begin{prop}\label{majodistance1}
Pour $k \in \Nat$ on a
\begin{equation*}
\delta_k \ioe \frac{1}{q_kq_{k+1}}
\quad\textrm{ et }\quad
\delta_k \soe\begin{cases}
  \frac{1/2}{q_{k+1}q_{k+2}} \quad &\text{si $a_{k+1}=1$}\\
  \frac{1/2}{q_kq_{k+1}} \quad &\text{si $a_{k+1}\soe 2$}.
\end{cases}
\end{equation*}
\end{prop}

\dem

On a d'abord, d'apr\`es \eqref{encabeta},
\begin{equation*}
\delta_k \ioe \frac{\beta_k}{q_k} \ioe \frac{1}{q_kq_{k+1}}.
\end{equation*}
Ensuite, si $a_{k+1}=1$, toujours d'apr\`es \eqref{encabeta},
\begin{equation*}
\delta_k = \frac{\beta_{k+1}}{q_{k+1}}  \soe \frac{1/2}{q_{k+1}q_{k+2}},
\end{equation*}
et si  $a_{k+1} \soe 2$,
\begin{equation*}
\delta_k \soe  \min \bigl (\frac{1}{q_{k}(q_k+q_{k+1})},\frac{1}{q_{k+1}(q_k+q_{k-1})}\bigr ) \soe \frac{1/2}{q_kq_{k+1}}.\qedhere
\end{equation*}

\medskip

La proposition suivante sera utilis\'ee au \S \ref{Le cas $k>K$} lors de la preuve de la proposition \ref{cas_ksupK}.

\begin{prop}\label{majodistance2}
Pour $k \in \Nat$ on a
$$
\frac{\log \big ((2\delta_k)^{-1}\big)}{q_k} \ioe 2\frac{\log q_{k+1}}{q_k} + 2\frac{\log
q_{k+2}}{q_{k+1}}.
$$
\end{prop}
\dem

Si $a_{k+1} \soe 2$ on a, d'apr\`es la proposition \ref{majodistance1},
\begin{align*}
 \frac{\log \big ((2\delta_k)^{-1}\big )}{q_k}  \ioe \frac{\log(q_kq_{k+1})}{q_k}
\ioe 2  \frac{\log q_{k+1}}{q_k}.
\end{align*}
Si $a_{k+1} =1$, on a
$$
\frac{\log \big ((2\delta_k)^{-1}\big )}{q_k}  \ioe \frac{\log(q_{k+1}q_{k+2})}{q_k}.
$$
Mais dans ce cas on a aussi $q_{k+1}\ioe 2q_k$, d'o\`u
$$
\frac{\log q_{k+2}}{q_k} \ioe 2 \frac{\log q_{k+2}}{q_{k+1}},
$$
ce qui démontre l'in\'egalit\'e annonc\'ee.\fin

\subsection{Profondeur d'un nombre rationnel}\label{profondeur_rationnel}

Soit $r$ un nombre rationnel, $0<r<1$, mis sous forme irr\'eductible $p/q$. Il peut s'\'ecrire d'une et une seule fa\^ion sous la forme
$$
r=[0; b_1,\dots,b_k]
$$
avec $k \in \Nat^*$, $b_i \in \Nat^*$, $1\ioe i\ioe k$, et $b_k \soe 2$. Nous dirons alors que $r$ est de profondeur\footnote{On trouvera dans \cite{MR1252067} une liste de probl\`emes arithm\'etiques sur la profondeur des nombres rationnels.} $k$. Par convention, $0$ et $1$ sont de profondeur $0$.

\'Ecrivons $[0; b_1,\dots,b_{k-1}]$ sous forme r\'eduite $p_{k-1}/q_{k-1}$ (si $k=1$, on a $p_0=0$ et $q_0=1$). Le nombre rationnel $r$ est une extr\'emit\'e de deux cellules de profondeur $k$ (qui sont donc contig\'ees) :
$$
\cgot \quad \text{d'extr\'emit\'es} \quad [0; b_1,\dots,b_{k-1},b_k-1]=\frac{p-p_{k-1}}{q-q_{k-1}}  \quad \text{et} \quad r \, ;
$$
$$
\cgot' \quad \text{d'extr\'emit\'es} \quad r \quad \text{et} \quad [0; b_1,\dots,,b_{k-1},b_k+1]=\frac{p+p_{k-1}}{q+q_{k-1}} .
$$

Comme $b_k \soe 2$, on a $q_{k-1} \ioe q/2$. La longueur de $\cgot$ est
$$
\frac{1}{(q-q_{k-1})q}>\frac{1}{q^2}.
$$
tandis que la longueur de $\cgot'$ est
$$
\frac{1}{q(q+q_{k-1})}\soe\frac{2}{3q^2}.
$$
En particulier, on a l'inclusion
$$
]r-\frac{2}{3q^2},r+\frac{2}{3q^2}[ \, \setminus \, \{r\} \subset \cgot \cup \cgot'.
$$

\subsection{Profondeur et \'epaisseur d'un segment inclus dans $]0,1[$}\label{t85}

Soit $I=[a,b]$ un segment inclus dans $]0,1[$, de longueur $h=b-a>0$ et d'extr\'emit\'es $a,b$ irrationnelles. Il existe un unique entier naturel $K$ tel que $I$ soit inclus dans une cellule de profondeur $K$, mais dans aucune cellule de profondeur $K+1$. Nous dirons alors que $I$ est de profondeur $K$. Dans le cas o\`u $x=(a+b)/2$ est aussi irrationnel, le nombre $K$ est d\'efini par l'encadrement
\begin{equation}
  \label{t84}
\delta_{K+1}(x) <h/2 < \delta_K(x).
\end{equation}

Signalons une petite subtilit\'e. Il est exact que la longueur d'un segment $I$ tend uniform\'ement vers $0$ quand sa profondeur $K$ tend vers l'infini : on a d'apr\`es \eqref{t83} et \eqref{mino_qk}
$$
|I| \ioe \frac{1}{F_{K+1}F_{K+2}}.
$$
En revanche, $|I|$ peut tendre vers $0$ alors que $K$ garde une valeur
constante, $0$ par exemple : il suffit que $1/2$ appartienne \`a $I$. Cela \'etant,
dans le cas o\`u $x \in X$ est fix\'e et $I=[x-h/2,x+h/2]$ avec $x\pm h/2\in X$,
on a
\begin{equation}\label{t100}
 \lim_{h\to 0} K(x,h)=\infty
\end{equation}
puisque
$$
\frac{1/2}{q_{K+2}(x)q_{K+3}(x)}<h
$$
d'apr\`es \eqref{t84} et la proposition \ref{majodistance1}.

Nous d\'efinissons \'egalement l'\'epaisseur de $I$ comme le nombre de cellules de profondeur $K+1$ qui ont une intersection non vide avec $I$, o\`u $K$ est la profondeur de $I$. L'\'epaisseur de $I$ est un nombre entier sup\'erieur ou \'egal \`a $2$.

\subsection{Convention}\label{t20}

Nous redéfinissons maintenant, et pour toute la suite, les fonctions $a_k$, $p_k$, $q_k$, $\alpha_k$, $\beta_k$ en les prolongeant par continuit\'e (sans changer de notation) sur chaque cellule de profondeur $k$ (et donc \'egalement sur chaque cellule de profondeur sup\'erieure \`a $k$). Ainsi, en posant $\cgot =\cgot(b_1,\dots,b_k)$, on aura pour $x \in \cgot$ :
\begin{align}
a_k(x)&=b_k \notag\\
\frac{p_k(x)}{q_k(x)}&=\frac{p_k}{q_k}=[0;b_1,\dots,b_k] \notag\\
\alpha_k(x)&=\frac{q_kx-p_k}{-q_{k-1}x+p_{k-1}} \label{alpha_cellule}\\
\beta_k(x)&=(-1)^{k-1}(p_k-xq_k).\label{beta_cellule}
\end{align}

En particulier, $\alpha_k$ est dérivable sur $\cgot$ et y v\'erifie
\begin{equation}\label{t65}
\alpha_k'=(-1)^{k}\beta_{k-1}^{-2}=(-1)^{k}(q_k+\alpha_kq_{k-1})^2.
\end{equation}

On a aussi la relation valable sur $\cgot$
\begin{equation}\label{t70}
\beta_{k-1}=\frac{1}{q_k}\bigl (1-q_{k-1}\beta_k\bigr).
\end{equation}

\section{Estimations en moyenne de $\gamma_k$ sur un intervalle inclus dans une cellule de profondeur~$K$}\label{pmajorations_moyenne}

Soit $I$ un intervalle inclus dans une cellule de profondeur $K$. Notre objectif est ici de majorer l'int\'egrale $\int_I\gamma_k(t)dt$, o\`u $\gamma_k$ est définie par \eqref{defgamma}. Pour $j\ioe K$, nous noterons $p_j$, $q_j$ les valeurs constantes de ces fonctions sur $I$.

La discussion portera esssentiellement sur les tailles respectives de $k$ et $K$. Nous ferons fr\'equemment usage, sans toujours le mentionner explicitement, des encadrements
\eqref{encabeta} et \eqref{encaalpha}. Enfin nous supposerons toujours que la longueur $h$ de $I$ v\'erifie $h \ioe e^{-2}$.

\subsection{Le cas $k<K$}

\begin{prop}\label{t79}
Pour $k<K$ on a
\begin{align*}
\int_I\gamma_k(t)dt &\ioe h \frac{\log (2q_{k+1})}{q_k}\\
&\ioe  \frac{h\log(1/h)}{2q_k} +\frac{h}{q_k}.
\end{align*}
\end{prop}
\dem

La premi\`ere in\'egalit\'e r\'esulte de \eqref{enca_gamma}. Pour la seconde, on observe que $I$ est inclus dans une cellule de profondeur $k+1$ dont la longueur est
$$
\frac{1}{q_{k+1}(q_{k+1}+q_k)}<\frac{1}{q_{k+1}^2},
$$
donc $q_{k+1} <h^{-1/2}$ et
\begin{align*}
\log (2q_{k+1}) &\ioe \log (h^{-1/2}) +\log 2\\
&\ioe \demi \log(1/h) + 1,
\end{align*}
d'o\`u le r\'esultat.\fin

\medskip

Par ailleurs, lorsque $k<K$, nous disposons d'une majoration uniforme de la dériv\'ee de $\gamma_k$ dans une cellule de profondeur $K$, ce qui nous conduit au r\'esultat suivant.
\begin{prop}\label{cas_kinfK}
Soit $I$ un segment de longueur $h$ inclus dans une cellule $\cgot$ de profondeur $K$ et $x$ un point quelconque de $I$. Pour $k<K$, on a
$$
 \int_I \lvert \gamma_k(t)-\gamma_k(x)\rvert dt \ioe \frac{3}{2} q_{k+1} h^2.
$$
\end{prop}
\dem

D'apr\`es les formules  \eqref{alpha_cellule} et \eqref{beta_cellule}, la fonction $\gamma_k$ est dérivable sur $\cgot$ pour $k\ioe K$, et \eqref{beta_cellule} et \eqref{t65} fournissent
\begin{equation}\label{derivgamma}
 \gamma_k'=(-1)^{k-1}q_{k-1}\log(1/\alpha_k)+\frac{(-1)^{k-1}}{\beta_k}.
\end{equation}
Nous en déduisons, pour $k<K$, $t\in I$,
\begin{equation}\label{majogamma'}
  \lvert \gamma_k'(t) \rvert \ioe q_{k-1}\log(a_{k+1}+1)+ q_{k+1}+q_k
\ioe q_k a_{k+1}+q_{k-1}+ 2q_{k+1} = 3 q_{k+1}.
\end{equation}
La conclusion découle alors de l'in\'egalit\'e des accroissements finis et de l'in\'egalit\'e
\begin{equation*}
\int_I \lvert t-x \rvert dt \ioe\frac{h^2}{2}. \qedhere
\end{equation*}


\subsection{Le cas $k=K$}\label{Le cas $k=K$}

Dans ce paragraphe et le suivant, nous donnons deux estimations, l'une comparant $\int_I \gamma_K(t) dt$ \`a $h\log 1/h$, l'autre \`a $h$. Pour simplifier les formulations et les d\'emonstrations qui suivent, nous supposons dans ces \S\S \ref{Le cas $k=K$} et \ref{Le cas $k>K$} que $I$ est un segment dont les extr\'emit\'es et le milieu sont irrationnels. Nous noterons respectivement $h$, $x$, $K$ et $E$ la longueur, le milieu, la profondeur et l'\'epaisseur de $I$.

Pour la premi\`ere estimation, un lemme \'el\'ementaire sera utile.

\begin{lem}\label{accroissements-t-logt}
  On a
  $$
  \int_I \log(1/t) dt
  \ioe h \log (1/h)+h.
  $$
\end{lem}
\dem

Comme $t \mapsto \log (1/t)$ est positive et d\'ecroissante sur $]0,1]$, on a
\begin{align*}
\int_I \log(1/t) dt &\ioe \int_0^{h} \log(1/t) dt\\
&= h \log(1/h)+h.\qedhere
\end{align*}

\begin{prop}\label{t71}
On a
$$
\int_I \gamma_K(t) dt \ioe \frac{h\log(1/h) }{q_K}([E>2]+[E=2]/2)
+\frac{2h}{q_K}.
$$
\end{prop}
\dem

En effet,
\begin{align*}
\int_I \gamma_K(t)dt&=
\int_I \beta_{K-1}(t) \log \beta_{K-1}(t) dt
+\int_I \beta_{K-1}(t) \log 1/\beta_{K}(t) dt\\
& \ioe \int_I \beta_{K-1}(t) \log 1/\beta_{K}(t) dt \expli{car $\beta_{K-1}(t)\ioe 1$}\\
&=\frac{1}{q_K^2}\int_{\beta_K(I)}
 (1-uq_{K-1})\log(1/u) du,
\end{align*} 
o\`u l'on a pos\'e $u=\beta_{K}(t)$  et utilis\'e
\eqref{t70}. On a donc
\begin{align*}
\int_I \gamma_K(t)dt &\ioe \frac{1}{q_K^2}\int_{\beta_K(I)}
 \log(1/u) du\\
&\ioe \frac{1}{q_K^2} |\beta_K(I)|\log 1/ |\beta_K(I)| + \frac{|\beta_K(I)|}{q_K^2}\expli{d'apr\`es le lemme \ref{accroissements-t-logt}}\\
&= \frac{1}{q_K^2}(q_Kh)\log((q_Kh)^{-1}) + \frac{q_Kh}{q_K^2}\\
&=\frac{h}{q_K} \log(1/h) + \frac{(1-\log q_K)h}{q_K},
\end{align*}
ce qui d\'emontre la proposition dans le cas o\`u $E>2$.

Par ailleurs, si $E=2$, les nombres $x\pm h/2$ se situent chacun dans l'une
de deux cellules adjacentes, disons $\cgot(b_1,\ldots,b_K,b_{K+1}-1)$ et
$\cgot(b_1,\ldots,b_K,b_{K+1})$, o\`u $b_{K+1}\soe 2$. Posons
$$
[0;b_1,\ldots,b_K,b_{K+1}]=\frac{p}{q} \quad ((p,q)=1),
$$
et $q$ est d'ailleurs la valeur maximale de la fonction $q_{K+1}$ sur $I$.

La longueur de $I$ est inf\'erieure \`a la somme des longueurs de ces deux cellules :

$$h \ioe \Big| \frac{p+p_K}{q+q_K}
- \frac{p-p_K}{q-q_K}\Big| =\frac{2}{(q+q_K)(q-q_K)}.$$
Comme $b_{K+1}\soe 2$, on a $q=b_{K+1}q_K+q_{K-1}\soe 2 q_K$, ce qui implique
\begin{equation}\label{majoration-ecart-xy}
h\ioe\frac{8/3}{q^2},
\end{equation}
et donc
$$
\log (2q) \ioe \demi\log (1/h)+\log (2\sqrt{8/3})\ioe \demi\log (1/h)+2.
$$

Enfin, nous avons
\begin{align*}
\int_I \gamma_K(t)dt&\ioe \frac{1}{q_K}\int_I \log \bigl (2q_{K+1}(t)\bigr)dt \expli{d'apr\`es \eqref{enca_gamma}}\\
&\ioe \frac{h}{q_K}\log (2q)  \\
&\ioe \frac{h}{2q_K}\log(1/h) +\frac{2h}{q_K}.\qedhere
\end{align*}

\begin{prop}\label{cas_kegalK}
On a
$$
\int_I \gamma_K(t) dt \ioe 8h \gamma_K(x)+\frac{h}{q_K} (6 \log q_K +4).
$$
\end{prop}
\dem

D'apr\`es  \eqref{encabeta}, on a
\begin{align*}
\int_I \gamma_K(t) dt &=\int_I  \beta_{K-1}(t) \log \bigl (1/\alpha_K(t)\bigr ) dt\\
&=\int_I  \beta_{K-1}(t)^3 \cdot  \beta_{K-1}(t)^{-2} \log \bigl (1/\alpha_K(t)\bigr ) dt\\
&\ioe q_K^{-3}\int_I  \beta_{K-1}(t)^{-2} \log \bigl (1/\alpha_K(t)\bigr ) dt.\end{align*}
Comme $\alpha_K'(t)= (-1)^K/\beta_{K-1}^2(t)$, il suit
\begin{equation}
\int_I \gamma_K(t) dt \ioe q_K^{-3}\int_{\alpha_K(I)}\log (1/u)du.
\end{equation}
L'intervalle $\alpha_K(I)$ a pour extr\'emit\'es $\alpha_K(x- h/2)$ et $\alpha_K(x+h/2)$. Or, d'apr\`es l'in\'egalit\'e des accroissements finis,
\begin{equation*}
  |\alpha_K(x\pm h/2)-\alpha_K(x)|\ioe(h/2) \max_{\xi\in I} \beta_{K-1}(\xi)^{-2} \ioe 2q_K^2h.
\end{equation*}
Supposons d'abord $4q_K^2h \ioe \alpha_K(x)$. Nous en déduisons
$$
\frac{1}{u} \ioe \frac{2}{\alpha_K(x)} \quad (u \in \alpha_K(I)),
$$
donc,
\begin{align*}
\int_I \gamma_K(t) dt &\ioe q_K^{-3}\log \bigl (2/\alpha_K(x)\bigr )   \int_{\alpha_K(I)}du\ioe q_K^{-3}\log \bigl (2/\alpha_K(x)\bigr ) \cdot 4q_K^2h\\
& \ioe 8h \beta_{K-1}(x)\log \bigl (1/\alpha_K(x)\bigr )+\frac{4h}{q_K}\log 2 \\
& \ioe 8h \gamma_K(x)+\frac{3h}{q_K}.
\end{align*}

Supposons ensuite $4q_K^2h \soe \alpha_K(x)$. Nous avons dans ce cas, d'apr\`es la proposition~\ref{majo_log_alphak},
\begin{align*}
  \int_I \gamma_K(t) dt &\ioe \max_{t\in I} |\beta_{K-1}(t)| \,e h\log (1/h)
 \ioe  \frac{e}{q_K}h\log \bigl (4q_K^2/\alpha_K(x) \bigr )\\
&=  \frac{e}{q_K}h\log \bigl (1/\alpha_K(x) \bigr )+ \frac{h}{q_K}(2e\log 2 + 2e\log q_K)\\& \ioe 6h \gamma_K(x)+ \frac{h}{q_K}(6\log q_K+4),
\end{align*}
d'o\`u la conclusion.\fin

\subsection{Le cas $k>K$}\label{Le cas $k>K$}

Nos calculs n\'ecessiteront une majoration \'el\'ementaire.
\medskip

\begin{lem}\label{somme_inverse_cubes}
Soit $m$ et $n$ des nombres entiers tels que $1 \ioe m <n$. On a alors
$$
\sum_{m \ioe \ell \ioe n}\frac{1}{\ell^3}  \ioe 3 \frac{n-m}{m^2n}.
$$
\end{lem}
\dem

On a
\begin{align*}
\sum_{m \ioe \ell \ioe n}\frac{1}{\ell^3}  &\ioe \frac{1}{m^3}+\int_m^n \frac{dt}{t^3}
 = \frac{1}{m^3} + \demi \Bigl ( \frac{1}{m^2}- \frac{1}{n^2} \Bigr )\\
&=  \frac{1}{m^3} + \frac{n^2-m^2}{2m^2n^2}
 \ioe \frac{1}{m^3} + \frac{n-m}{m^2n}.
\end{align*}
Maintenant, comme $1 \ioe m \ioe n-1$, on a $
  m(n-m) \soe n-1\soe n/2$,
donc
$$
\frac{1}{m^3}\ioe 2 \frac{n-m}{m^2n},
$$
d'o\`u le r\'esultat.\fin

\bigskip

Nous sommes maintenant en mesure de traiter le cas $k>K$.

\begin{prop}\label{t72}
Pour $k \in \Nat$ tel que $k> K$ on a
$$
\int_I\gamma_{k}(t)dt \ioe [E=2] \frac{6}{q_{K+1}(x)F_{k-K}} h \log 1/h + [E>2]\frac{ 72}{q_KF_{k-K}}h.
$$
\end{prop}
\dem

Soit $\cgot$ la cellule de profondeur $K$ qui contient $I$. Comme $\cgot$ est l'intervalle d'extr\'emit\'es
$$
\frac{p_K}{q_K} \quad \text{et} \quad \frac{p_K+p_{K-1}}{q_K+q_{K-1}},
$$
nous pouvons \'ecrire
\begin{align*}
  x+(-1)^Kh/2 &=  \frac{up_K+p_{K-1}}{uq_K+q_{K-1}}\\
  x+(-1)^{K-1}h/2 &=  \frac{vp_K+p_{K-1}}{vq_K+q_{K-1}},
\end{align*}
avec
$1\ioe u<1/\alpha_K(x)<v<+\infty$. Posons $m=\lfloor u \rfloor$ et $n=\lfloor v \rfloor$, de sorte que $1 \ioe m \ioe a_{K+1}(x) \ioe n$ et $E=n-m+1$, $n>m$.

On a donc, en utilisant l'identit\'e \eqref{pkqk1},
\begin{align}
  h &= \left \lvert  \frac{up_K+p_{K-1}}{uq_K+q_{K-1}} -\frac{vp_K+p_{K-1}}{vq_K+q_{K-1}} \right \rvert\notag\\
&= \frac{v-u}{(uq_K+q_{K-1})(vq_K+q_{K-1})} \notag\\
& \soe \frac{v-u}{q_K^2(u+1)(v+1)} \notag\\
& \soe \frac{v-u}{6q_K^2mn},\label{mino_h}
\end{align}
o\`u la derni\`ere minoration découle de $u+1\ioe 3m $ et $v+1\ioe 2n$.

Rappelons que les cellules de profondeur $K+1$ incluses dans $\cgot$ sont les intervalles
$$
\cgot_\ell = \Bigl \{\frac{sp_{K}+p_{K-1}}{sq_{K}+q_{K-1}}, \, \ell<s<\ell+1 \Bigr \} \quad (\ell \soe 1).
$$

On a
\begin{align*}
x+(-1)^Kh/2 &\in \overline{\cgot_m},\\
  x+(-1)^{K-1}h/2 & \in \overline{\cgot_n}.
\end{align*}

D'autre part, $I$ est inclus dans la r\'eunion des $ \overline{\cgot_\ell}$, $m \ioe \ell \ioe n$, donc
\begin{equation}\label{t101} 
\int_I \gamma_{k}(t)dt \ioe \sum_{n\ioe \ell\ioe m} \int_{\cgot_\ell} \beta_{k-1}(t)\log\big(1/\alpha_k(t)\big)dt
.
\end{equation} 
Distinguons alors deux cas. 

\textit{Premier cas} : $E=2$.

On a donc $n=m+1$. Comme $I \subset \overline {\cgot_m} \cup  \overline
{\cgot_{m+1}}$, on a
$$
q_{K+1}(t) =
\begin{cases}
  mq_K+q_{K-1}\\
\text{ou}\\
 (m+1)q_K+q_{K-1}
\end{cases}
\quad (t \in I, \, t \not \in \partial\cgot_m )
$$
En particulier, pour $t \in I$, $t \not \in \partial\cgot_m$, on a
\begin{equation}\label{t73}
  \frac{q_{K+1}(t)}{q_{K+1}(x)} \soe \frac{mq_K+q_{K-1}}{(m+1)q_K+q_{K-1}}\soe
\demi.
\end{equation}
Par cons\'equent, d'apr\`es la majoration \eqref{majobetakj} et la proposition
\ref{majo_log_alphak},
\begin{align}
\int_I \gamma_{k}(t)dt &\ioe  \frac{2}{q_{K+1}(x)F_{k-K}}\int_I
\log\big(1/\alpha_k(t)\big) dt \ioe \frac{2}{q_{K+1}(x)F_{k-K}}
eh\log(1/h).\label{t74}
 \end{align}

\textit{Deuxi\`eme cas} : $E \soe 3$.

On a donc $v-u \soe 1$, d'o\`u
\begin{equation}\label{t102} 
n-m \ioe v-(u-1)\ioe 2(v-u).
\end{equation}
D'apr\`es la majoration \eqref{majobetakj} appliqu\'ee avec $i=K$ et $j=k-K-1$ dans
\eqref{t101}, 
\begin{align*}
\int_I \gamma_{k}(t)dt &\ioe \frac{1}{F_{k-K} }\sum_{n\ioe \ell\ioe m} \int_{\cgot_\ell}
\frac{\log\Bigl (1/\alpha_{k-K-1}\bigl (\alpha_{K+1}(t)\bigr )\Bigr )}{q_{K+1}(t)} d t \\
& \ioe \frac{1}{q_KF_{k-K} }\sum_{n\ioe \ell\ioe m}\frac{1}{\ell} \int_{\cgot_\ell}
\log \Bigl (1/\alpha_{k-K-1}\bigl (\alpha_{K+1}(t)\bigr )\Bigr ) d t,
\end{align*}
o\`u la derni\`ere in\'egalit\'e provient du fait que pour $t\in \cgot_\ell$, $q_{K+1}(t)=\ell q_K+q_{K-1}$.
Nous effectuons dans chaque $ \int_{\cgot_\ell}$ le changement de variable $w=\alpha_{K+1}(t)$. Notant que d'apr\`es \eqref{t65}
$$dt = (-1)^{K+1} \frac{ d w}{(q_{K+1}+wq_K)^2},$$ nous obtenons,
\begin{align*}
\int_I \gamma_{k}(t)dt &\ioe \frac{1}{q_KF_{k-K} }\sum_{m\ioe \ell\ioe n}\frac{1}{\ell} \int_0^{1}
\log\big(1/\alpha_{k-K-1}(w)\big) \frac{d w}{(q_{K+1}+wq_K)^2} \\
&\ioe   \frac{1}{q_K^3F_{k-K}}\int_0^1 \log\big(1/\alpha_{k-K-1}(w)\big) dw\,
\sum_{m\ioe \ell \ioe n} \frac{1}{\ell^3}\\
& \ioe   \frac{2}{q_K^3F_{k-K}}\int_0^1 \log\big(1/\alpha_{k-K-1}(w)\big) \frac{dw}{1+w}\,
\sum_{m\ioe \ell \ioe n} \frac{1}{\ell^3}.
\end{align*}
En utilisant le lemme \ref{somme_inverse_cubes} et l'invariance de la mesure $\mu$ par $\alpha_{k-K-1}$, il vient
\begin{equation*}
\int_I \gamma_{k}(t)dt\ioe \frac{6(n-m)}{q_K^3F_{k-K}m^2n}
\int_0^1 \log(1/w)dw =
 \frac{6(n-m)}{q_K^3F_{k-K}m^2n},
\end{equation*}
et, d'apr\`es \eqref{mino_h} et \eqref{t102},
\begin{equation} \label{t76}
\int_I \gamma_{k}(t)dt  \ioe 12\frac{v-u}{q_{K}^3F_{k-K}m^2n}
 \ioe 72 \frac{h}{q_KF_{k-K}m} \ioe 72 \frac{h}{q_KF_{k-K}}.
\end{equation} 

La conclusion d\'ecoule de  \eqref{t74} et \eqref{t76}. \fin

\begin{prop}\label{cas_ksupK}
Pour $k \in \Nat$ tel que $k> K$ on a
\begin{align*}
\int_I&\gamma_{k}(t)dt\\& \ioe \frac{h}{F_{k-K}}\Bigl (
\frac{72[E>2]}{q_K}+12\frac{[E=2]\log q_{K+2}(x)}{q_{K+1}(x)}+12\frac{[E=2]\log
q_{K+3}(x)}{q_{K+2}(x)}\Bigr ).
\end{align*}
\end{prop}
\dem

On a
\begin{align*}
 \frac{\log(1/h)}{q_{K+1}(x)} &\ioe \frac{\log
\bigl((2\delta_{K+1}(x))^{-1}\bigr)}{q_{K+1}(x)} \\
&\ioe 2\frac{\log q_{K+2}(x)}{q_{K+1}(x)}+2\frac{\log q_{K+3}(x)}{q_{K+2}(x)}
\end{align*}
d'apr\`es la proposition \ref{majodistance2}. Le r\'esultat d\'ecoule de cette in\'egalit\'e et de la proposition \ref{t72}.\fin

\section{Comportement de $\Psi$ au voisinage d'un point rationnel}\label{t77}

Nous commen\c{c}ons par traiter les cas des points $0$ et $1$.

\begin{lem}\label{t90}
On a
$$
\int_0^x \Phi\bigl (\alpha(t)\bigr)dt =O(x) \quad (0\ioe x\ioe 1).
$$
\end{lem}
\dem

Pour $0<x<1$, on a
\begin{align*}
\int_0^x \Phi\bigl (\alpha(t)\bigr)dt &\ioe \int_0^{1/a_1(x)}\Phi\bigl (\alpha(t)\bigr)dt \expli{car $x \ioe 1/\lfloor 1/x\rfloor =1/a_1(x)$}\\
&=\sum_{n\soe a_1(x)}\int_{1/(n+1)}^{1/n}\Phi\bigl (\tfrac{1}{t}-n\bigr)dt
\\
&=\sum_{n\soe a_1(x)}\int_0^1\Phi(u)\frac{du}{(n+u)^2}\\
&=\int_0^1\Phi(u)\sum_{n\soe a_1(x)}\frac{1}{(n+u)^2}du\\
&\ioe \frac{2}{a_1(x)}\int_0^1\Phi(u)du =O(x).\qedhere
\end{align*}

\begin{lem}\label{ratio_zero}
Pour $0<x \ioe 1$, on a
\begin{equation}\label{cas_zero}
\Psi(x)=x\log (1/x) +x+O(x^2),
\end{equation}
et
\begin{equation}\label{cas_un}
\Psi(1)-\Psi(1-x)=x\log (1/x) +x+O\bigl(x^2\log(2/x)\bigr)
 \end{equation}
\end{lem}

\dem

Pour $0<x \ioe 1$, on a, d'apr\`es l'\'equation fonctionnelle de $\Phi$,
\begin{align*}
\Psi(x)&=\int_0^x \Phi  (t) dt
= \int_0^x \Bigl (t\Phi\bigl( \alpha(t)\bigr)+\log (1/t)\Bigr )dt \\
&= \int_0^x \log (1/t)dt  + O \Bigl(x\int_0^x\Phi\bigr( \alpha(t)\bigl )dt\Bigr)=x\log (1/x) +x+O(x^2),
\end{align*}
d'apr\`es le lemme \ref{t90}.

Pour \'etablir \eqref{cas_un}, il suffit de traiter le cas $0<x<1/2$. On a
\begin{align*}
\Psi(1)-\Psi(1-x)&=\int_{1-x}^1 \Phi  (t) dt
= \int_{1-x}^1 \Bigl (t\Phi\bigl( \alpha(t)\bigr)+\log (1/t)\Bigr )dt \\
&= \int_{1-x}^1 \Bigl (t\Phi\bigl( (1-t)/t\bigr)+\log (1/t)\Bigr )dt.
\end{align*}
En effectuant le changement de variable $u=(1-t)/t$ puis en employant \eqref{cas_zero}, nous obtenons
\begin{align*}
\Psi(1)-\Psi(1-x)&= \int_0^{x/(1-x)} \Bigl ( \frac{\Phi(u)}{1+u} +\log(1+u)\Bigr )\frac{du}{(1+u)^2}\\
&= \int_0^{x/(1-x)} \Bigl ( \Phi(u)+O\bigl (u\Phi(u)\bigr ) +O(u)\Bigr )du\\
&=\frac{x}{1-x}\log\frac{1-x}{x} +\frac{x}{1-x}+O(x^2\log(1/x))\\
&=x\log (1/x) +x+O\bigl (x^2\log (1/x)\bigr ).\qedhere
\end{align*}

Nous sommes maintenant en mesure d'\'etablir le point $(iii)$ du th\'eor\`eme \ref{t68}.
\begin{prop}\label{cas_rationnel}
On a uniform\'ement pour $r$ rationnel et $h$ r\'eel tels que $0<r<1$, $q^2|h|<2/3$ (o\`u $q$ est le dénominateur de l'\'ecriture irr\'eductible de $r$) :
\begin{equation}\label{t66}
\Psi(r+h)-\Psi(r)=\frac{h}{q} \log (1/|h|) +\Big(\frac{1}{q}-2 \frac{\log
q}{q}+ \Phi(r)\Big)h +O\bigl(qh^2 \log ((q^2|h|)^{-1}) \bigr).
\end{equation}

\end{prop}

\dem

 Soit $K \in \Nat^*$ la profondeur du nombre rationnel $r$. Comme $0<|h|<2/3q^2$, on sait d'apr\`es le \S \ref{profondeur_rationnel} que l'intervalle ouvert d'extr\'emit\'es $r$ et $r+h$ est inclus dans l'une des deux cellules de profondeur $K$ :

$\bullet$ $\cgot$ d'extr\'emit\'es $\frac{p-p_{K-1}}{q-q_{K-1}}$ et $r$ ;

$\bullet$ $\cgot'$ d'extr\'emit\'es $r$ et $\frac{p+p_{K-1}}{q+q_{K-1}}$.

Les trois nombres $\frac{p-p_{K-1}}{q-q_{K-1}}$, $r$ et $\frac{p+p_{K-1}}{q+q_{K-1}}$ se succ\'edent dans cet ordre si $K$ est pair, dans l'ordre inverse si $K$ est impair. Par cons\'equent,
\begin{equation*}
  (-1)^Kh>0 \Leftrightarrow r+h\in \cgot'.
\end{equation*}

On a
\begin{align*}
 \Psi(r+h)-\Psi(r) &= \int_r^{r+h}\Phi(t)dt\\
&=\sum_{k<K} \int_r^{r+h}\gamma_k(t)dt+\int_r^{r+h}\beta_{K-1}(t)\Phi\bigl(\alpha_K(t)\bigr)dt,
\end{align*}
d'apr\`es l'identit\'e \eqref{t61}. Or,
\begin{align*}
\left \lvert \int_r^{r+h}\sum_{k<K}\gamma_k(t)dt-h\Phi(r) \right \rvert
&=\left \lvert \sum_{k<K}  \int_r^{r+h}\bigl (\gamma_k(t)-\gamma_k(r)\bigr )dt \right \rvert\\
&\ioe \frac{3}{2}h^2 \sum_{k<K} q_{k+1} \ioe \frac 92 qh^2,
\end{align*}
d'apr\`es la proposition \ref{cas_kinfK} et l'in\'egalit\'e \eqref{somme_qk}.

Il reste donc \`a \'etablir
\begin{equation}\label{queue_rationnel}
\int_r^{r+h}\beta_{K-1}(t)\Phi\bigl(\alpha_K(t)\bigr)dt=\frac{1}{q} h\log
((q^2|h|)^{-1}) +\frac{h}{q}+O\bigl(qh^2 \log ((q^2|h|)^{-1}) \bigr).
\end{equation}

La fonction $\alpha_K$ est définie sur chacune des cellules $\cgot$ et $\cgot'$. Suivant que $r+h$ appartienne \`a $\cgot$ ou \`a $\cgot'$, nous prolongeons $\alpha_K$ (par continuit\'e \`a droite ou \`a gauche suivant la parit\'e de $K$) respectivement \`a $\cgot \cup \{r\}$ ou $\cgot' \cup \{r\}$, et dans l'int\'egrale du premier membre de \eqref{queue_rationnel} nous effectuons le changement de variables $u=\alpha_K(t)$. Distinguons les deux cas.

$\bullet$ Premier cas : $(-1)^Kh >0$. Dans ce cas, l'intervalle d'extr\'emit\'es $r$ et $r+h$ est inclus dans $\overline{\cgot'}$. On a $q_K=q$, $p_K=p$ et

\begin{equation*}
  u=\alpha_K(t)=\frac{q_Kt-p_K}{-q_{K-1}t
  +p_{K-1}},
\end{equation*}
de sorte que, compte tenu de \eqref{t65},
\begin{align*}
  \int_r^{r+h}\beta_{K-1}(t)\Phi\bigl(\alpha_K(t)\bigr)dt
  &=(-1)^K\int_{\alpha_K(r)}^{\alpha_K(r+h)}
  \Phi(u)\frac{du}{(q_K+q_{K-1}u)^3}.
\end{align*}

On a $\alpha_K(r)=0$ et
\begin{align*}
  \alpha_K(r+h)
  &=\frac{q^2h}{-q_{K-1}p-q_{K-1}qh+p_{K-1}q}\\
&=\frac{q^2|h|}{1-|h|qq_{K-1}}.
\end{align*}
Nous posons donc
\begin{equation}\label{def_x'}
  x'=\frac{q^2|h|}{1-|h|qq_{K-1}}.
\end{equation}
Comme $q_{K-1}\ioe q/2$, on a $0<x'<1$ et
\begin{align*}
\int_r^{r+h}&\beta_{K-1}(t)\Phi\bigl(\alpha_K(t)\bigr)dt
\\&=\frac{(-1)^K}{q^3}
\int_0^{x'}
\Phi(u) \frac{du}{\Big(1+u\frac{q_{K-1}}{q}\Big)^3}\\
&=\frac{(-1)^K}{q^3}\int_0^{x'}
\Phi(u)\bigl (1+O(u)\bigr )du\\
&=\frac{(-1)^K}{q^3}x'\log(1/x')+\frac{(-1)^K}{q^3}x'+O\bigl(q^{-3}x'^2\log (2/x')\bigr)
\end{align*}
o\`u la derni\`ere \'egalit\'e r\'esulte de \eqref{cas_zero} et du lemme \ref{t90}.
\medskip

$\bullet$ Deuxi\`eme cas : $(-1)^K h<0$.
 Dans ce cas, l'intervalle d'extr\'emit\'es $r$ et $r+h$ est inclus dans $\overline{\cgot}$. Dans cette cellule, on a
\begin{equation*}
  u=\alpha_K(t)=\frac{(q-q_{K-1})t-p+p_{K-1}}{-q_{K-1}t
  +p_{K-1}},
\end{equation*}
Par suite, $\alpha_K(r)=1$ et
\begin{align*}
  \alpha_K(r+h)
  &= \frac{(q-q_{K-1})\frac{p}{q}-(p-p_{K-1})+h(q-q_{K-1})}{-q_{K-1}\frac{p}{q}+p_{K-1}
  -q_{K-1}h}\\
  &=1+\frac{hq^2}{(-1)^K-qq_{K-1}h}=1-x,
\end{align*}
avec cette fois
\begin{equation}\label{def_x}
  x=\frac{q^2|h|}{1+|h|qq_{K-1}}\, (\ioe 2/3).
\end{equation}

Nous avons donc, d'apr\`es \eqref{cas_un},
\begin{align*}
\int_r^{r+h}&\beta_{K-1}(t)\Phi\bigl(\alpha_K(t)\bigr)dt
\\&=(-1)^{K-1}
\int_{1-x}^1 \Phi(u) \frac{du}{\big(q+(u-1)q_{K-1}\big)^3}\\
&=\frac{(-1)^{K-1}}{q^3}
\int_{1-x}^1 \Phi(u) du + O\Big(q^{-3}\int_{1-x}^1 \Phi(u)(1-u)du\Big)\\
&=\frac{(-1)^{K-1}}{q^3}
\int_{1-x}^1 \Phi(u) du +O\bigl(q^{-3}x^2\log (1/x)\bigr)\\
&=\frac{(-1)^{K-1}}{q^3}x\log(1/x)+\frac{(-1)^{K-1}}{q^3}x+O\bigl(q^{-3}x^2\log (1/x)\bigr).
\end{align*}

Nous avons donc dans les deux cas, avec $y=x$ ou $x'$, 
\begin{equation*}
\int_r^{r+h}\beta_{K-1}(t)\Phi\bigl(\alpha_K(t)\bigr)dt=
\frac{\sgn (h)}{q^3}y\log(1/y)+\frac{\sgn (h)}{q^3}y+O\bigl(q^{-3}y^2\log
(2/y)\bigr)
\end{equation*}
Comme $q_{K-1}<q/2$ et $|h|<2/3q^2$, nous avons
\begin{equation*}
  \log(1/y)=\log((q^2|h|)^{-1}) +O(q^2h),
\end{equation*}
et
\begin{equation*}
  y=q^2|h|+O(h^2 q^4),
\end{equation*}
ce qui entra\^ine bien \eqref{queue_rationnel}.\fin

\bigskip
\goodbreak

\section{Comportement de $\Psi$ au voisinage d'un nombre de Cremer}\label{etude_Psi}

Nous commen\c{c}ons par d\'emontrer le point $(i)$ du th\'eor\`eme \ref{t68}, qui concerne tous les irrationnels, qu'ils soient de Cremer ou de Brjuno.

\begin{prop}\label{irrationnel_gene}
Soit $x\in X$. On a
\begin{equation*}
 \Psi(x+h/2)-\Psi(x-h/2)=o(h\log(1/h)) \quad(h\to0).
\end{equation*}
\end{prop}

\dem

Soit $x\in X$ et $K \in \Nat$. Soit $0<h<\min\big(e^{-2},2\delta_K(x)\big)$ de sorte que l'intervalle $I=]x-h/2,x+h/2[$ soit inclus dans la cellule de profondeur $K$ contenant $x$.

Nous avons, d'apr\`es les majorations \eqref{enca_gamma} et \eqref{majobetakj},
\begin{equation*}
\int_I \Phi  (t) dt
\ioe h\sum_{k<K}\frac{\log (2q_{k+1})}{q_k}+ \sum_{k\soe K} \frac{1}{q_K F_{k-K+1}}\int_I \log\big(\alpha_k(t)\big)dt,
\end{equation*}
avec $q_k=q_k(x)$ pour $k\ioe K$.
En employant la proposition \ref{majo_log_alphak} et l'in\'egalit\'e \eqref{t82}, nous obtenons
\begin{equation*}
\int_I \Phi  (t) dt
 \ioe 4 h\log (2q_K) + \frac{10}{q_K} h \log(1/h).
\end{equation*}
Par cons\'equent,
$$
\limsup_{h \vers 0} \frac{\lvert \Psi(x+h/2)-\Psi(x-h/2)\rvert}{h \log (1/h)} \ioe \frac{10}{q_K}.
$$
Comme $K$ est arbitraire, cela d\'emontre le r\'esultat.\fin
\bigskip

Nous d\'emontrons maintenant le point $(ii)$ du th\'eor\`eme \ref{t68} lorsque $x$ est un nombre de Cremer.
\begin{prop}\label{cremer}
 Soit $x$ un nombre de Cremer. On a
\begin{equation*}
\lim_{h\to 0} \frac{\Psi(x+h)-\Psi(x)}{h}=+\infty.
\end{equation*}
\end{prop}

\dem

 Comme la fonction $\gamma_k$ est positive, nous avons en toute g\'en\'eralit\'e pour $h>0$,
\begin{align*}
\frac{\Psi(x+h)-\Psi(x)}{h}&=\frac{1}{h} \int_x^{x+h} \sum_{k \soe 0} \gamma_k(t)dt
\soe \frac{1}{h} \int_x^{x+h} \sum_{k \ioe K} \gamma_k(t)dt\\
&=\sum_{k \ioe K} \frac{ \int_x^{x+h} \gamma_k(t)dt }{h}.
\end{align*}
Si $x$ est irrationnel, chacune des fonctions $\gamma_k$ est continue au point $x$, ce qui entra\^ine
\begin{equation*}
\liminf_{h\vers 0} \frac{\Psi(x+h)-\Psi(x)}{h} \soe \sum_{k\ioe K} \gamma_k(x) \quad(x\in X,K\in\Nat).
\end{equation*}
Si $x$ est un point de Cremer, le membre de droite tend vers l'infini avec $K$, ce qui implique le r\'esultat voulu.
\fin
\bigskip

\section{Démonstration du th\'eor\`eme \ref{thprincipal} et comportement de
$\Psi$ au voisinage d'un nombre de Brjuno }\label{pth_principal}

Les propositions \ref{cremer} et \ref{cas_rationnel} impliquent que tout point
de Lebesgue de $\Phi$ est n\'e\-ces\-sai\-re\-ment un nombre de Brjuno. Il reste
donc \`a \'etablir que
\begin{equation}
  \label{t86}
  \frac{1}{h}\int_{x-h/2}^{x+h/2} \lvert\Phi(t)-\Phi(x) \rvert dt\vers 0 \quad
(h\vers 0),
\end{equation}
quand $x$ est un nombre de Brjuno. Cela impliquera directement le 
point (ii) du th\'eor\`eme \ref{t68} lorsque $x_0$ est un nombre de Brjuno.

Par densit\'e de $\Real \setminus \Rat$ dans
$\Real$, on peut supposer que les nombres $x\pm h/2$ sont irrationnels. Posons
alors $I=[x-h/2,x+h/2]$ et notons $K=K(x,h)$ la profondeur de ce segment.
Nous avons
\begin{equation*}
\begin{split}
\int_I \lvert\Phi(t)-\Phi(x) \rvert dt
&\ioe \sum_{k<K} \int_I \lvert \gamma_k(t)-\gamma_k(x) \rvert dt
+ \int_I \gamma_K(t) dt
 \\
&\quad+\sum_{k> K} \int_I \gamma_k(t) dt+h\sum_{k\soe K} \gamma_k(x).
\end{split}
\end{equation*}

Majorons chacun des quatre termes du second membre. En notant désormais
$q_j=q_j(x)$ pour tout $j$, on a
\begin{align*}
\frac{1}{h}\sum_{k<K} \int_I \lvert \gamma_k(t)-\gamma_k(x) \rvert dt &\ioe \frac{3h}{2} \sum_{k<K} q_{k+1} \expli{d'apr\`es la proposition \ref{cas_kinfK}}\\
& \ioe \frac 92 q_K h \quad \text{\footnotesize (d'apr\`es \eqref{somme_qk})}\\
&\ioe \frac{9}{q_{K+1}},
\end{align*}
la derni\`ere in\'egalit\'e r\'esultant de la proposition \ref{majodistance1} puisque
$h/2<\delta_K(x)$) ; puis,
\begin{align*}
\frac{1}{h}  \int_I \gamma_K(t) dt &\ioe  8 \gamma_K(x)+\frac{6 \log q_K +4}{q_K} \quad \text{\footnotesize (d'apr\`es la proposition \ref{cas_kegalK})}\\
&\ioe \frac{14 \log q_{K+1} +10}{q_K} \quad\text{\footnotesize (d'apr\`es \eqref{enca_gamma})} \, ;
\end{align*}
ensuite, d'apr\`es la proposition \ref{cas_ksupK},
\begin{align*}
\frac{1}{h} \sum_{k> K} \int_I \gamma_k(t) dt &\ioe \Bigl (
\frac{72}{q_K}+12\frac{\log q_{K+2}}{q_{K+1}}+12\frac{\log
q_{K+3}}{q_{K+2}}\Bigr ) \sum_{k> K}\frac{1}{F_{k-K}} \\
& \ioe  \frac{242}{q_K}+41\frac{\log q_{K+2}}{q_{K+1}}+41\frac{\log q_{K+3}}{q_{K+2}} ;
\end{align*}
et enfin
\begin{align*}
\sum_{k\soe K} \gamma_k(x) \ioe \sum_{k\soe K} \frac{\log(2q_{k+1})}{q_{k}}
\quad \text{\footnotesize (d'apr\`es \eqref{enca_gamma})}.
\end{align*}

Finalement, on obtient
\begin{equation}\label{majolongue}
\frac{1}{h}\int_I \lvert\Phi(t)-\Phi(x) \rvert dt \ioe \sum_{k\soe K} \frac{42\log q_{k+1}+253}{q_{k}}
\end{equation}
Puisque $x$ est un nombre de Brjuno, le membre de droite de \eqref{majolongue}
tend vers $0$ lorsque $K$ tend vers l'infini, et donc lorsque $h$ tend vers $0$
d'apr\`es \eqref{t100}, ce qui ach\'eve la preuve.

\section{D\'emonstration du th\'eor\`eme \ref{t69}}\label{module-continuite}

Le lemme \ref{ratio_zero} du \S \ref{t77} prouve que
$$
\omega(h)\soe h\log (1/h) +O(h).
$$
Il reste donc \`a \'etablir
\begin{equation}\label{majoration-module-continuite}
\omega(h) \ioe h\log (1/h) +O(h).
\end{equation}

Lorsque $h>e^{-2}$, il s'agit de montrer que $\omega(h)=O(1)$, ce qui est vrai car $\Phi \in L^1(0,1)$. Il suffit donc de montrer que si $I$ est un segment de longueur $h\ioe e^{-2}$ inclus dans $[0,1]$ on a
\begin{equation}\label{module-continuite-a-montrer}
\int_I\Phi(t)dt
\ioe h\log(1/h)+O(h).
\end{equation}
En fait, par densit\'e de $\Real \setminus \Rat$ dans $\Real$, on peut \'egalement supposer que les extr\'emit\'es et le milieu $x$ de $I$ sont irrationnels. Soit alors $K$ la profondeur de $I$ et $E$ son \'epaisseur.

Notons $\cgot=\cgot(a_1,\ldots,a_K)$ la cellule de profondeur $K$ qui contient $I$ et, pour $j\ioe K$, désignons par $p_j$ et $q_j$ les valeurs constantes de ces fonctions dans $\cgot$. Comme $I$ n'est pas inclus dans une cellule de profondeur $K+1$, il existe $a_{K+1}\soe 2$ tel que le nombre rationnel
$$
r=[0;a_1,\ldots,a_K,a_{K+1}]=\frac{p}{q} \quad ((p,q)=1)
$$
appartienne \`a $I$. Observons que $p$ et $q$ sont les valeurs
respectives des fonctions $p_{K+1}$ et $q_{K+1}$ dans la cellule
$\cgot(a_1,\ldots,a_K,a_{K+1})$. D'autre part, si $E=2$, le nombre $x$ se situe soit dans la cellule $\cgot(a_1,\ldots,a_{K},a_{K+1})$, soit dans la cellule $\cgot(a_1,\ldots,a_{K},a_{K+1}-1)$. Par cons\'equent
\begin{equation}\label{t89}
q_{K+1}(x)=q \;\text{\rm ou} \; (q-q_K),\; \text{et donc}\; q_{K+1}(x)\soe q/2.
\end{equation}

\medskip

Nous distinguons alors deux cas : $q\ioe 100$ et $q>100$, o\`u la valeur $100$ est choisie pour obtenir simplement l'inégalité cherchée dans le second cas (pour le premier cas, n'importe quelle valeur convient).

\textit{Premier cas : $q\ioe 100$.}

\'Ecrivons $I=[r -\alpha h,r+\beta h]$ avec $0<\alpha,\beta <1$ et $\alpha
+\beta=1$. On peut supposer $h<1/15000$, et donc $q^2h<2/3$. D'apr\`es la
proposition \ref{cas_rationnel}, on a alors
\begin{align*}
\int_I\Phi(t)dt&=  \bigl (\Psi(r+\beta h)-\Psi(r)\bigr )+\bigl (\Psi(r)-\Psi( r -\alpha h)\bigr )\\
&= \frac{\beta h}{q} \log ((\beta h)^{-1})
+\frac{\alpha h}{q}\log ((\alpha h)^{-1}) +O(h)\\
&= \frac{h \log(1/h)}{q} + \frac{h}{q} \big(\alpha \log(1/\alpha)
+ \beta \log(1/\beta) \big) + O(h)\\
&\ioe h\log 1/h+O(h),
\end{align*}
car la fonction $t\mapsto t\log(1/t)$ est born\'ee sur l'intervalle $[0,1]$.

\goodbreak
\medskip

\textit{Deuxi\`eme cas : $q> 100$.}

En notant $m$ le plus grand entier inf\'erieur ou \'egal \`a $K$ tel que
$q_m\ioe 100$, nous effectuons la décomposition
$$
\int_I\Phi(t)dt= \sum_{k<m} \int_I \gamma_k(t) dt
+ \int_I \gamma_m(t)dt+ \sum_{k>m} \int_I \gamma_k(t) dt.$$
Nous avons, d'apr\`es \eqref{enca_gamma},
$$
\sum_{k<m} \int_I \gamma_k(t) dt
\ioe  h \sum_{k<m} \frac{\log (2q_{k+1})}{q_k}\ioe h \log (2q_m) \sum_{k\soe 1} \frac{1}{F_k}
\ioe  3,36h \log 200,$$
et donc
$$
\int_I\Phi(t)dt\ioe \int_I \gamma_m(t)dt+ \sum_{k>m} \int_I \gamma_k(t) dt+O(h).$$
Il reste  \`a \'etablir la majoration
\begin{equation}
  \label{t78}
\int_I \gamma_m(t)dt+ \sum_{k>m} \int_I \gamma_k(t) dt \ioe h\log(1/h) +O(h).
\end{equation}

Pour ce faire, nous distinguons deux sous-cas.
\smallskip

\textit{Premier sous-cas : $m<K$.}

La fonction $q_{m+1}$ est alors constante sur $I$ et, par maximalit\'e de $m$, on a
$q_{m+1}>100$. Nous avons alors directement avec \eqref{majobetakj} et la proposition
\ref{majo_log_alphak}
\begin{align*}
  \sum_{k>m} \int_I \gamma_k(t) dt
&\ioe \frac{1}{q_{m+1}}\sum_{k>m}\frac{1}{F_{k-m}}
\int_I \log(1/\alpha_{k}(t)) dt\\
&\ioe \frac{3,36 e h \log(1/h)}{100} \ioe \frac{h}{4}\log(1/h).
\end{align*}
Par ailleurs la proposition \ref{t79} nous donne
$$
\int_I \gamma_m(t)dt\ioe \frac{h\log(1/h)}{2} +h.
$$

On obtient donc
$$ \int_I\Phi(t)dt\ioe \frac{3}{4}h\log(1/h) +O(h).$$

\smallskip

\textit{Deuxi\`eme sous-cas : $m=K$.}

Nous appliquons alors les propositions \ref{t71} et \ref{t72} :
\begin{align*}
\int_I \gamma_K(t)dt+& \sum_{k>K} \int_I \gamma_k(t) dt\\ & \ioe h\log
1/h\Bigl (  \frac{[E>2]+[E=2]/2}{q_K}  +
\frac{6[E=2]}{q_{K+1}(x)}\sum_{k>K}\frac{1}{F_{k-K}}\Bigr)\\&\quad +
h\Bigl (\frac{1}{q_K}+ \frac{ 72[E>2]}{q_K}\sum_{k>K}\frac{1}{F_{k-K}} \Bigr).
\end{align*}

Comme
\begin{align*}
\frac{[E>2]+[E=2]/2}{q_K}&  +
\frac{6[E=2]}{q_{K+1}(x)}\sum_{k>K}\frac{1}{F_{k-K}} \\&\ioe [E>2]+[E=2]\Bigl
(\demi+\frac{12}{q}3,36\Bigr )\expli{d'apr\`es \eqref{t89}}\\
&\ioe 1\expli{car $q>100$},
\end{align*}
nous avons bien \'etabli \eqref{module-continuite-a-montrer} dans tous les cas.
\fin

\section{Questions}\label{t87}

Il serait d'abord int\'eressant de déterminer le comportement asymptotique du module de continuit\'e $L^1$ de la fonction de Brjuno,
$$
\omega_1(h)=\sup_{0<h'\ioe h}\int_0^1|\Phi(t+h')-\Phi(t)|dt
$$
(o\`u $\Phi$ est prolong\'ee par $1$-p\'eriodicit\'e sur $\Real$), quand $h$ tend vers $0$.

\smallskip

D'autre part, dans \cite{Marmi_Moussa_Yoccoz-2}, Marmi, Moussa et Yoccoz prolongent la fonction de Brjuno au demi-plan $\Im z >0$. Ils étudient en détail les propriétés de la \emph{fonction de Brjuno complexe}, et notamment son comportement quand $\Im z \vers 0$. Il serait intéressant de relier nos résultats et les leurs, en particulier le point $(iii)$ du théorème \ref{t68} avec le \S 5.2.9 de \cite{Marmi_Moussa_Yoccoz-2}. Comme nous avons conçu le présent travail dans le cadre strict de l'analyse réelle, nous laissons ces questions ouvertes.

\medskip

\begin{multicols}{2}
\footnotesize

\noindent BALAZARD, Michel\\
Institut de Math\'ematiques de Luminy, UMR 6206\\
CNRS, Aix Marseille Université\\
Campus de Luminy, Case 907\\
13288 Marseille Cedex 9\\
FRANCE\\
Adresse \'electronique : \\ \texttt{balazard@iml.univ-mrs.fr}

\smallskip

\noindent MARTIN, Bruno\\
Laboratoire de Math\'ematiques Pures et Appliqu\'ees\\
CNRS, Universit\'e du Littoral C\^ote d'Opale\\
50 rue F. Buisson, BP 599\\
62228 Calais Cedex\\
FRANCE\\
Adresse \'electronique :  \\ \texttt{martin@lmpa.univ-littoral.fr}
\end{multicols}

\end{document}